\documentclass[dvips, 10pt, a4paper]{article}
\usepackage{latexsym, amsmath}
\usepackage{amssymb}

\begin{document}
\sloppy
\centerline{\Large{A geometric inequality for warped product semi-slant}}

\smallskip

\centerline{\Large{submanifolds of nearly cosymplectic manifolds}}

\bigskip

\centerline{{Siraj Uddin, Abdulqader Mustafa, Bernardine R. Wong and Cenap Ozel}}

\bigskip

\begin{abstract}
Recently, we have shown that there do not exist the warped product semi-slant submanifolds of cosymplectic manifolds [10]. As nearly cosymplectic structure generalizes cosymplectic ones same as nearly Kaehler generalizes Kaehler structure in almost Hermitian setting. It is interesting that the warped product semi-slant submanifolds exist in nearly cosymplectic case while in case of cosymplectic do not exist. In the beginning, we prove some preparatory results and finally we obtain an inequality such as $\|h\|^2 \geq 4q\csc^2\theta\{1+\frac{1}{9}\cos^2\theta\}\|\nabla \ln f\|^2$ in terms of intrinsic and extrinsic invariants. The equality case is also considered.\\


\end{abstract}

\section{Introduction}
To study the manifolds with negative curvature Bishop and O'Neill introduced the idea of warped products [2]. Afterwards, this idea was used to model the standard space time, especially in the neighborhood of massive stars and black holes [11]. However, warped product spaces were developed in Riemannian manifolds enormously [8, 9]. The geometry of warped product submanifolds is intensified after B.Y. Chen's work on warped product CR-submanifolds of Kaehler manifolds [8]. Motivated by Chen's papers, many geometers studied warped product submanifolds in almost Hermitian as well as almost contact metric manifolds [1, 9, 13, 15].

\parindent=8mm
On the other hand, almost contact manifolds with Killing structures tensors were defined in [3] as nearly cosymplectic manifolds.  Later on, Blair and Showers [5] studied nearly cosymplectic structure $(\phi, \xi, \eta, g)$ on a Riemannian manifold  $\bar M$ with $\eta$ closed from the topological viewpoint. An almost contact metric structure $(\phi, \xi, \eta, g)$ satisfying  $(\bar\nabla_X\phi)X=0$ is called a nearly cosymplectic structure. If we consider $S^5$ as a totally geodesic hypersurface of $S^6$; then it is known that $S^5$ has a non cosymplectic nearly cosymplectic structure. It was shown that the normal nearly cosymplectic manifolds are cosymplectic (see [4]).

\parindent=8mm Next, the slant submanifolds of an almost contact metric manifold were defined and studied by J.L. Cabrerizo et.al [7]. The notion of semi-slant submanifolds of almost Hermitian manifolds was introduced by N. Papaghuic [12]. In fact, semi-slant submanifolds in almost Hermitian manifolds are defined on the line of CR-submanifolds. In the setting of almost contact metric manifolds, semi-slant submanifolds are defined and investigated by Cabrerizo et. al [6].

\parindent=8mm Recently, we studied warped product semi-slant submanifolds of cosymplectic manifolds [10]. We have seen that there do not exist warped product semi-slant submanifolds in cosymplectic manifolds. As the nearly cosymplectic structure is generalizes the cosymplectic ones in the same sense as nearly Kaehler generalizes Kaehler structure. Therefore, the geometric study of warped product semi-slant submanifolds in nearly cosymplectic is interesting. Such type of warped products exist in nearly cosymplectic case while in case of cosymplectic do not exist. In this paper, we prove that the warped product semi-slant submanifolds of the type $N_\theta\times{_{f}N_T}$ do not exist in a nearly cosymplectic manifold $\bar M$. However, we obtain some interesting results on the existence of the warped product submanifolds of the type $N_T\times{_{f}N_\theta}$ of a nearly cosymplectic manifold $\bar M$, where $N_T$ and $N_\theta$ are the invariant and proper slant submanifolds of $\bar M$, respectively. We also establish a general sharp inequality for squared norm of the second fundamental form in terms of the warping function and the slant angle for the warped product semi-slant submanifolds in the form $N_T\times{_{f}N_\theta}$ in a nearly cosymplectic manifold $\bar M$.

\section{Preliminaries}
A $(2n+1)-$dimensional $C^\infty$ manifold $\bar M$ is said to have an {\it{almost contact structure}} if there exist on $\bar M$ a tensor field $\phi$ of type $(1,1)$, a vector field $\xi$ and a $1-$form $\eta$ satisfying [5]
$$\phi^2=-I+\eta\otimes\xi,~~\phi\xi=0,~~\eta\circ\phi=0,~~\eta(\xi) = 1.\eqno(2.1)$$
There always exists a Riemannian metric $g$ on an almost contact manifold $\bar M$ satisfying the following compatibility condition
$$\eta(X)=g(X, \xi),~~~g(\phi X, \phi Y)=g(X, Y)-\eta(X)\eta(Y)\eqno(2.2)$$
where $X$ and $Y$ are vector fields on $\bar M$ [5].

\parindent=8mm
An almost contact structure $(\phi, \xi, \eta)$ is said to be {\it normal} if the almost complex structure $J$ on the product manifold $\bar M\times\mathbb R$ given by
$$J(X, f\frac{d}{dt})=(\phi X-f\xi,~\eta(X)\frac{d}{dt}),$$
where $f$ is a $C^\infty-$function on $\bar M\times\mathbb R$ has no torsion i.e., $J$ is integrable, the condition
for normality in terms of $\phi,~\xi$ and $\eta$ is $[\phi, \phi]+2d\eta\otimes\xi=0$ on $\bar M$, where $[\phi, \phi]$ is the Nijenhuis tensor of $\phi$. Finally the {\it{fundamental 2-form}} $\Phi$ is defined by $\Phi(X, Y)=g(X, \phi Y)$.

\parindent=8mm
An almost contact metric structure $(\phi, \xi, \eta, g)$ is said to be {\it{cosymplectic}}, if
it is normal and both $\Phi$ and $\eta$ are closed [5]. The structure is said to be {\it{nearly cosymplectic}} if $\phi$ is Killing, i.e., if
$$(\bar\nabla_X\phi)Y+(\bar\nabla_Y\phi)X=0,\eqno(2.3)$$
for any $X,~Y\in T\bar M$, where $T\bar M$ is the tangent bundle of $\bar M$ and $\bar\nabla$ denotes the Riemannian connection of the metric $g$. Equation (2.3) is equivalent to $(\bar\nabla_X\phi)X=0$, for each $X\in T\bar M$. The structure is said to be {\it{closely cosymplectic}} if $\phi$ is Killing and $\eta$ is closed. It is well known that an almost contact metric manifold is {\it{cosymplectic}} if and only if $\bar\nabla\phi$ vanishes identically, i.e., $(\bar\nabla_X\phi)Y=0$ and $\bar\nabla_X\xi=0$.\\

\noindent{\bf{Proposition 2.1 [5].}} {\it{On a nearly cosymplectic manifold the vector field $\xi$ is
Killing.}}

\parindent=8mm
From the above proposition, we have $g(\bar\nabla_X\xi, X)=0$, for any vector field $X$ tangent to $\bar M$, where $\bar M$ is a nearly cosymplectic manifold.

\parindent=8mm
Let $M$ be submanifold of an almost contact metric manifold $\bar M$ with induced metric $g$ and let $\nabla$ and $\nabla^{\perp}$ be the induced connections on the tangent bundle $TM$ and the normal bundle $T^{\perp}M$ of $M$, respectively. Denote by ${\cal{F}}(M)$ the algebra of smooth functions on $M$ and by $\Gamma(TM)$ the ${\cal{F}}(M)$-module of smooth sections of $TM$ over $M$. Then the Gauss and Weingarten formulas are given by
$$\bar\nabla_XY=\nabla_XY+h(X, Y)\eqno(2.4)$$
$$\bar\nabla_XN=-A_NX+\nabla^{\perp}_XN,\eqno(2.5)$$
for each $X,~Y\in\Gamma(TM)$ and $N\in \Gamma(T^\perp M)$, where $h$ and $A_N$ are the second
fundamental form and the shape operator (corresponding to the normal vector
field $N$) respectively for the immersion of $M$ into $\bar M$. They are related as
$$g(h(X, Y), N)=g(A_NX, Y),\eqno(2.6)$$
where $g$ denotes the Riemannian metric on $\bar M$ as well as induced on $M$.

\parindent=8mm
For any $p\in M$, let $\{ e_1,\cdots, e_m,\cdots, e_{2n+1}\}$ be an orthonormal frame for the tangent space $T_p\bar M$, such that $e_1,\cdots, e_m$ are tangent to $M$ at $p$. We denote by $H$ the mean curvature vector, that is
$$H(p)=\frac {1}{m}\sum_{i=1}^mh(e_i, e_i).$$
Also, we set
$$h_{ij}^r = g( h(e_i, e_j), e_r), ~~~ i, j \in \{1,\cdots, m\}, ~r\in \{m+1,\cdots, 2n+1\}$$
and
$$||h||^2 = \sum _{i, j = 1}^m g( h(e_i, e_j), h(e_i, e_j)).$$

\parindent=8mm
For any $X\in \Gamma(TM)$, we write
$$\phi X=PX+FX,\eqno(2.7)$$
where $PX$ is the tangential component and $FX$ is the normal component of $\phi X$.

\parindent=8mm
Similarly for any $N\in \Gamma(T^\perp M)$, we write
$$\phi N=BN+CN,\eqno(2.8)$$
where $BN$ is the tangential component and $CN$ is the normal component of $\phi N$. The covariant derivative of the tensor field $\phi$ is defined as
$$(\bar\nabla_X\phi)Y=\bar\nabla_X\phi Y-\phi\bar\nabla_XY,\eqno(2.9)$$
for all $X,~Y\in \Gamma(T\bar M)$.

\parindent=8mm
Now, denote by ${\cal P}_XY$ and ${\cal Q}_XY$ the tangential and normal
parts of $(\bar\nabla_X\phi)Y$, i.e.,
$$(\bar\nabla_X\phi)Y={\cal P}_XY+{\cal Q}_XY\eqno(2.10)$$
for all $X, Y\in \Gamma(TM)$. Making use of (2.7)-(2.10) and the Gauss and Weingarten formulae, the following
equations may easily be obtained
$${\cal P}_XY=(\bar\nabla_XP)Y-A_{FY}X-Bh(X,Y)\eqno(2.11)$$
$${\cal Q}_XY=(\bar\nabla_XF)Y+h(X, PY)-Ch(X, Y)\eqno(2.12)$$
where the covariant derivative of $P$ and $F$ are defined by
$$(\bar\nabla_X P)Y=\nabla_XPY-P\nabla_XY\eqno(2.13)$$
$$(\bar\nabla_XF)Y=\nabla^\perp_XFY-F\nabla_XY\eqno(2.14)$$
for all $X,Y\in \Gamma(TM)$.

\parindent=8mm
Similarly, for any $X\in\Gamma(TM)$ and $N\in \Gamma(T^\perp M)$, denoting the tangential and
normal parts of $(\bar\nabla_X\phi)N$ by ${\cal P}_XN$ and ${\cal Q}_XN$
respectively, we obtain
$${\cal P}_XN=(\bar\nabla_XB)N+PA_NX-A_{CN}X \eqno(2.15)$$
$${\cal Q}_XN=(\bar\nabla_XC)N+h(BN, X)+FA_NX \eqno(2.16)$$
where the covariant derivative of $B$ and $C$ are defined
by
$$(\bar\nabla_XB)N=\nabla_XBN-B\nabla^\perp_XN\eqno(2.17)$$
$$(\bar\nabla_XC)N=\nabla^\perp_XCN-C\nabla^\perp_XN.\eqno(2.18)$$

\parindent=8mm
It is straightforward to verify  the following properties of ${\cal P}$ and ${\cal Q}$, which we enlist here for later use:
$$(p_1)~~(i)~~~~{\cal P}_{X + Y}W={\cal P}_XW + {\cal P}_YW,
~~~~~~(ii)~~~{\cal Q}_{X + Y}W={\cal Q}_XW + {\cal Q}_YW,~~~$$
$$~(p_2)~~~(i)~~~~{\cal P}_X(Y+W)={\cal P}_XY + {\cal P}_XW,~~~~(ii)~~~~{\cal Q}_X(Y + W)={\cal Q}_XY + {\cal Q}_XW,$$
$$(p_3)~~(i)~~~~g({\cal P}_XY,~ W)=- g(Y, {\cal P}_XW),~~~(ii)~~~~g({\cal Q}_XY,~ N)=-g(Y, {\cal P}_XN),$$
for all $X, Y, W\in\Gamma(TM)$ and $N\in\Gamma(T^\perp M)$.

\parindent=8mm
On a submanifold $M$ of a nearly cosymplectic  manifold, by equations (2.3) and (2.10) we have
$$(a)~~{\cal P}_XY + {\cal P}_YX=0,~~(b)~~{\cal Q}_XY + {\cal Q}_YX= 0\eqno(2.19)$$
for any $X,~Y\in\Gamma(TM)$.

\parindent=8mm
A submanifold $M$ of an almost contact metric manifold $\bar M$ is said to be {\it{invariant}} if $F$ is identically zero, that is, $\phi X\in \Gamma(TM)$ and {\it{anti-invariant}} if $P$ is identically zero, that is, $\phi X\in \Gamma(T^\perp M)$, for any $X\in \Gamma(TM)$.

\parindent=8mm
We shall always consider $\xi$ to be tangent to the submanifold $M$. There is another class of submanifolds that is called the slant submanifold. For each non zero vector $X$ tangent to $M$ at $x$, such that $X$ is not proportional to $\xi_x$, we denote by $0\leq\theta (X) \leq\pi /2$, the angle between $\phi X$ and $T_x M$ is called the slant angle. If the slant angle $\theta (X)$ is constant for all $X\in T_x M - \langle\xi_x\rangle$ and $\ x\in M$, then M is said to be {\it{slant submanifold}} [7]. Obviously if $\theta = 0,~M$ is invariant and if $\theta$ = $\pi /2,~M$ is an anti-invariant submanifold. A slant submanifold is said to be {\it{proper slant}} if it is neither invariant nor anti-invariant submanifold.

\parindent=8mm
We recall the following result for a slant submanifold.\\

\noindent
{\bf {Theorem 2.1 [7]}} {\it Let M be a submanifold of an almost contact metric manifold $\bar M$, such that $\xi$ is tangent to $M$. Then $M$ is slant if and only if there exists a constant  $\lambda\in [0, 1]$ such that
$$P^2=\lambda (- I + \eta \otimes \xi)\eqno (2.20)$$
Furthermore, if $\theta$ is slant angle of $M$, then $\lambda = \cos ^2\theta.$}

\parindent=8mm
Following relations are straightforward consequences of relation (2.20)
$$g(PX,PY)=\cos^2\theta ( g(X,Y)- \eta (Y)\eta (X))\eqno(2.21)$$
$$g(FX,FY)=\sin^2\theta ( g(X,Y)-\eta (Y)\eta (X))\eqno(2.22)$$
for all $X,Y\in \Gamma(TM)$.

\parindent=8mm
A submanifold $M$ of an almost contact metric manifold $\bar M$ is said to be a {\it semi-slant} if there exist two orthogonal distributions
$D_1$ and $D_2$ satisfying:
\begin{enumerate}
\item[{(i)}] $TM={\cal D}_1\oplus {\cal D}_2\oplus \langle\xi\rangle$
\item[{(ii)}] ${\cal D}_1$ is an invariant  i.e., $\phi {\cal D}_1\subseteq {\cal D}_1$
\item[{(iii)}] ${\cal D}_2$ is a slant distribution with slant angle $\theta\neq \frac{\pi}{2}.$
\end{enumerate}

\parindent=8mm
A semi-slant submanifold $M$ of an almost contact  manifold $\bar M$ is {\it mixed geodesic} if
$$h(X, Z)=0\eqno(2.23)$$
for any $X\in {\cal D}_1$ and $Z\in {\cal D}_2$. Moreover, if $\mu$ is the $\phi-$invariant subspace of the normal bundle $T^\perp M$, then in case of semi-slant submanifold, the normal bundle $T^\perp M$ can be decomposed as
$$T^\perp M= F{\cal D}_2\oplus \mu.\eqno(2.24)$$

\section{Warped product semi-slant submanifolds}
Bishop and O'Neill [2] introduced the notion of warped product manifolds. These manifolds are natural generalizations of Riemannian product manifolds. They defined these manifolds as: Let $(N_1, g_1)$ and $(N_2, g_2)$ be two Riemannian manifolds and $f>0$ a differentiable function on $N_1$. Consider the product manifold $N_1\times N_2$ with its projections $\pi_1:N_1\times N_2\to N_1$ and $\pi_2:N_1\times N_2\to N_2$. Then the warped product of $N_1$ and $N_2$ denoted by $M=N_1\times{_f}N_2$ is a Riemannian manifold $N_1\times N_2$ equipped with the Riemannian structure such that
$$g(X, Y)=g_1({\pi_1}_\star X, {\pi_1}_\star Y)+(f\circ\pi)^2g_2({\pi_2}_\star X, {\pi_2}_\star Y)$$
for each $X, Y\in\Gamma(TM)$ and $\star$ is a symbol for the tangent map. Thus we have
$$g=g_1+f^2g_2.\eqno(3.1)$$
The function $f$ is called the {\it{warping function}} of the warped product [2, 11]. A warped product manifold $N_1\times{_f}N_2$ is said to be {\it{trivial}} if the warping function $f$ is constant. We recall the following general result for a warped product manifold for later use.\\

\noindent
{\bf {Lemma 3.1 [2]}} {\it Let $M = N_1\times_{f}N_2$ be a warped product manifold. Then}
\begin{enumerate}
\item[{(i)}] {\it$\nabla_XY\in TN_1$ is the lift of $\nabla_XY$ on $N_1$}
\item[{(ii)}] {\it $\nabla_XZ=\nabla_ZX=(X\ln f)Z$}
\item[{(iii)}] {\it $\nabla_ZW = \nabla^{N_2}_ZW - g(Z, W) \nabla\ln f$}
\end{enumerate}
{\it{for each X, Y$\in\Gamma(TN_1)$ and $Z, W\in \Gamma(TN_2)$ where $\nabla$ and $\nabla ^{N_2}$ denote the Levi-Civita connections on $M$ and $N_2$, respectively, and $\nabla\ln f$ is the gradient of $\ln f$.}}

\parindent=8mm Let $M$ be a Riemannian manifold of dimension $k$ with the inner product $g$ and $\{e_1,\cdots,e_k\}$ be an orthonormal frame on $M$. Then for a differentiable function $f$ on $M$, the gradient $\nabla f$ of a function $f$ on $M$ is defined by
$$g(\nabla f, U)=Uf,\eqno(3.2)$$
for any $U\in \Gamma(TM)$. As a consequence, we have
$$\|\nabla f\|^2=\sum\limits_{i=1}\limits^k (e_i(f))^2\eqno(3.3)$$
where $\nabla f$ is the gradient of the function $f$ on $M$.

\parindent=8mm
Now, we consider the warped product semi-slant submanifolds tangent to the structure vector field $\xi$ which are either in the form $M=N_T\times{_{f}N_\theta}$ or $M=N_\theta\times{_{f}N_T}$ in a nearly cosymplectic  manifold $\bar M$, where $N_T$ and $N_\theta$ are invariant and proper slant submanifolds of a nearly cosymplectic manifold $\bar M$, respectively. On a warped product submanifold $M=N_1\times{_{f}N_2}$ of a nearly cosymplectic manifold $\bar M$, we have the following result.\\

\noindent
{\bf{Theorem 3.1. [15]}} {\it A warped product submanifold $M=N_1\times{_{f} N_2}$ of a nearly cosymplectic manifold $\bar M$ is an usual Riemannian product if the structure vector field $\xi$ is tangent to $N_2$, where $N_1$ and $N_2$ are the Riemannian submanifolds of $\bar M$.}

\parindent=8mm
From the above theorem for the existence of warped products we always consider the structure vector field $\xi$ is tangent to the base. Now, we start with the warped product semi-slant submanifolds of the type $M=N_\theta\times{_{f}N_T}$ of a nearly cosymplectic manifold $\bar M$.\\

\noindent
{\bf{Theorem 3.2}} {\it There do not exist the warped product semi-slant submanifolds M=$N_\theta \times{_{f}N_T}$ of a nearly cosymplectic manifold $\bar M$, where $N_\theta$ and $N_T$ are proper slant and invariant submanifolds of $\bar M$, respectively.}\\

\noindent
{\it Proof.} The proof is similar to the semi-invariant case, which we have proved in [15].$~\blacksquare$

\parindent=8mm
Now, we discuss the other case and all results before the geometric inequality are preparatory and we can not skip all these. We shall use these preparatory results in our main theorem.\\

\noindent
{\bf{Lemma 3.2}} {\it Let $M=N_T\times{_{f}N_\theta}$ be a warped product semi-slant submanifold of a nearly cosymplectic manifold $\bar M$, where $N_T$ and $N_\theta$ are invariant and proper slant submanifolds of $\bar M$, respectively. Then}
\begin{enumerate}
\item[{(i)}] $\xi\ln f=0$
\item[{(ii)}] $g(h(X, Y), FZ)=0$
\item[{(iii)}] $g(h(PX, Z), FZ)=(X\ln f)\|Z\|^2$
\item[{(iv)}] $g(h(X, Z), FPZ)=-g(h(X, PZ), FZ)$
\item[{(v)}]  $g({\cal P}_XZ, PZ)= 2  g(h(X, Z), FPZ)$
\end{enumerate}

\noindent {\it{for any $X, Y\in\Gamma(TN_T)$ and $Z\in\Gamma(TN_\theta)$.}}\\

\noindent{\it Proof.} Let $M=N_T\times{_{f}N_\theta}$ be a warped product semi-slant submanifold of a nearly cosymplectic manifold $\bar M.$ We assume that the structure vector field $\xi$ is tangent to $N_T$, then for any $Z\in\Gamma(TN_\perp)$, we have
$$\bar\nabla_Z\xi=\nabla_Z\xi+h(Z, \xi).$$
Taking the inner product with $Z$ and using Proposition 2.1 and Lemma 3.1 $(ii)$, we obtain
$$(\xi\ln f)\|Z\|^2=0.$$
This means that either $M$ is invariant or $\xi\ln f=0$, which proves $(i)$. Now, we consider $X, Y\in\Gamma(TN_T)$ and $Z\in\Gamma(TN_\theta)$, the we have
$$g(h(X, Y), FZ)=g(\bar\nabla_XY, FZ)=-g(Y, \bar\nabla_XFZ).$$
Using (2.7) and then (2.9), we obtain
$$g(h(X, Y), FZ)=-g(Y, (\bar\nabla_X\phi)Z)-g(Y, \phi\bar\nabla_XZ)+g(Y, \bar\nabla_XPZ).$$
Then from (2.2), (2.4) and Lemma 3.1 (ii), the second and last terms of right hand side vanish identically and hence by (2.10), we derive
$$g(h(X, Y), FZ)=-g(Y, {\cal P}_XZ).$$
Thus, on using the property $p_3~(i)$, we get
$$g(h(X, Y), FZ)=g({\cal P}_XY, Z).$$
Hence, by skew-symmetry of ${\cal P}_XY$ and symmetry of $h(X, Y)$, we get the second part of the lemma. For the third part, consider for any $X\in\Gamma(TN_T)$ and $Z\in\Gamma(TN_\theta)$, we have
$$g(h(PX, Z), FZ)=g(\bar\nabla_{Z}PX, \phi Z-PZ)~~~~~~~~~~~~~~$$
$$~~~~~~~~~~~~~~~~~~~~~~~~~=-g(PX, \bar\nabla_{Z}\phi Z)-g(\bar\nabla_{Z}PX, PZ).$$
From (2.4), (2.9) and Lemma 3.1 $(ii)$, the above equation reduced to
$$g(h(PX, Z), FZ)=-g(PX, (\bar\nabla_{Z}\phi) Z)-g(PX, \phi\bar\nabla_{Z}Z)-(PX\ln f)g(Z, PZ).$$
On using the structure equation of nearly cosymplectic and the fact that $Z$ and $PZ$ are orthogonal vector fields, the first and last terms of right hand side are identically zero. Then from (2.2), we derive
$$g(h(PX, Z), FZ)=g(\phi^2X, \bar\nabla_{Z}Z).$$
Using (2.1), we get
$$g(h(PX, Z), FZ)=-g(X, \bar\nabla_{Z}Z)+\eta(X)g(\xi, \bar\nabla_{Z}Z).$$
By the property of Riemannian connection $\bar\nabla$, the above equation takes the form
$$g(h(PX, Z), FZ)=g(\bar\nabla_{Z}X, Z)-\eta(X)g(\bar\nabla_{Z}\xi, Z).$$
Then from (2.4), Proposition 2.1 and Lemma 3.1 $(ii)$, we obtain
$$g(h(PX, Z), FZ)=(X\ln f)\|Z\|^2,$$
which is third part of the lemma. For the other parts, consider
$$g(\nabla_{PZ}\phi X, Z)=g(\bar\nabla_{PZ}\phi X, Z)$$
for any $X\in\Gamma(TN_T)$ and $Z\in\Gamma(TN_\theta)$. Using the property of Riemannian connection $\bar\nabla$ and Lemma 3.1 $(ii)$, we get
$$(\phi X\ln f)g(Z, PZ)=-g(\phi X, \bar\nabla_{PZ}Z).$$
Using the fact that $Z$ and $PZ$ are orthogonal vector fields, the above equation reduces to
$$0=g(X, \phi\bar\nabla_{PZ}Z).$$
Then form (2.9), we derive
$$0=g(X, \bar\nabla_{PZ}\phi Z)-g(X, (\bar\nabla_{PZ}\phi)Z).$$
By (2.7) and (2.10), we obtain
$$0=g(X, \bar\nabla_{PZ}PZ)+g(X, \bar\nabla_{PZ}FZ)-g(X, {\cal{P}}_{PZ}Z).$$
Using (2.4), (2.5) and (2.19) $(a)$, we get
$$0=-g(\nabla_{PZ}X, PZ)-g(X, A_{FZ}PZ)+g(X, {\cal{P}}_{Z}PZ).$$
Then from the property $p_3~(i)$ and Lemma 3.1 $(ii)$, we obtain
$$0=-(X\ln f)g(PZ, PZ)-g(h(X, PZ), FZ)-g({\cal{P}}_{Z}X, PZ).$$
Again using (2.19) $(a)$, (2.21) and the fact that $\xi$ is tangent to $N_T$, we derive
$$g({\cal{P}}_{X}Z, PZ)=(X\ln f)\cos^2\theta\|Z\|^2+g(h(X, PZ), FZ).\eqno(3.4)$$
Interchanging $Z$ by $PZ$ and then using (2.20), (2.21) and the fact that $\xi$ is tangent to $N_T$, we obtain
$$-\cos^2\theta g({\cal{P}}_{X}PZ, Z)=(X\ln f)\cos^4\theta\|Z\|^2-\cos^2\theta g(h(X, Z), FPZ).$$
By the property $p_3~(i)$, the above equation will be
$$g({\cal{P}}_{X}Z, PZ)=(X\ln f)\cos^2\theta\|Z\|^2-g(h(X, Z), FPZ).\eqno(3.5)$$
Thus, the fourth part of the lemma follows from (3.4) and (3.5). Now, for the part $(v)$, we consider
$$g(h(X, PZ), FZ)=g(\bar\nabla_XPZ, FZ)$$
for any $X\in\Gamma(TN_T)$ and $Z\in\Gamma(TN_\theta)$. Using the property of Riemannian connection $\bar\nabla$ and then using (2.7), we have
$$g(h(X, PZ), FZ)=-g(PZ, \bar\nabla_X\phi Z)+g(PZ, \bar\nabla_XPZ).$$
Using (2.9), Lemma 3.1 $(ii)$, (2.21) and the fact that $\xi$ is tangent to $N_T$, we obtain
$$g(h(X, PZ), FZ)=-g(PZ, \phi\bar\nabla_XZ)-g(PZ, (\bar\nabla_X\phi)Z)+(X\ln f)\cos^2\theta\|Z\|^2.$$
Then from (2.2) and (2.10), we get
$$g(h(X, PZ), FZ)=g(\phi PZ, \bar\nabla_XZ)-g(PZ, {\cal{P}}_XZ)+(X\ln f)\cos^2\theta\|Z\|^2.$$
Using (2.4) and (2.7), we derive
$$g(h(X, PZ), FZ)=g(P^2Z, \nabla_XZ)+g(h(X, Z), FPZ)~~~~~~~~~~~~~~~$$
$$~~~~~~~~~~~~~~~-g(PZ, {\cal{P}}_XZ)+(X\ln f)\cos^2\theta\|Z\|^2.$$
Again, from the fact that $\xi$ is tangent to $N_T$ and using (2.20), the above equation reduces to
$$g({\cal{P}}_XZ, PZ)=g(h(X, Z), FPZ)-g(h(X, PZ), FZ).\eqno(3.6)$$
Thus, the fifth part of the lemma follows from (3.6) and the fourth part of this lemma. This proves the lemma completely.$~\blacksquare$\\

\noindent
{\bf{Theorem 3.3}} {\it Let $M=N_T\times{_{f}N_\theta}$ be a warped product semi-slant submanifold in a nearly cosymplectic manifold $\bar M$, where $N_T$ and $N_\theta$ are invariant and proper slant submanifolds of $\bar M$, respectively. Then,
$$g({\cal P}_XZ, PZ)=\frac{2}{3}(X\ln f)\cos^2\theta\|Z\|^2$$
for any $X\in\Gamma(TN_T)$ and $Z\in\Gamma(TN_\theta)$.}\\

\noindent
{\it Proof.} From (3.4) and (3.5), we have
$$2g({\cal{P}}_{X}Z, PZ)=2(X\ln f)\cos^2\theta\|Z\|^2+g(h(X, PZ), FZ)$$
$$-g(h(X, Z), FPZ).\eqno(3.7)$$
Then, by Lemma 3.2 $(iv)$, we obtain
$$2g({\cal{P}}_{X}Z, PZ)=2(X\ln f)\cos^2\theta\|Z\|^2-2g(h(X, Z), FPZ).\eqno(3.8)$$
Thus, from Lemma 3.2 $(v)$ and (3.8), we obtain the desire result. $~\blacksquare$

\parindent=8mm From the above theorem we have the following consequence.\\

\noindent
{\bf{Corollary 3.1}} {\it{ A warped product semi-slant submanifold $M=N_T\times{_{f}N_\theta}$ of a nearly cosymplectic manifold $\bar M$ is simply a Riemannian product of $N_T$ and $N_\theta$ if and only if ${\cal P}_XZ\in\Gamma(TN_T)$, for any X$\in\Gamma(TN_T)$ and Z$\in\Gamma(TN_\theta)$.}}\\

\noindent
{\bf{Theorem 3.4}} {\it Let $M=N_T\times{_{f}N_\theta}$ be a warped product semi-slant submanifold in a nearly cosymplectic manifold $\bar M$, where $N_T$ and $N_\theta$ are invariant and proper slant submanifolds of $\bar M$, respectively. Then,
$$g(h(X, Z), FPZ)=-g(h(X, PZ), FZ)=\frac{1}{3}(X\ln f)\cos^2\theta\|Z\|^2,$$
for any $X\in\Gamma(TN_T)$ and $Z\in\Gamma(TN_\theta)$.}\\

\noindent
{\it Proof.} The first equality is nothing but Lemma 3.2 (iv) and the second equality is directly followed by the equation (3.8) and Lemma 3.2 $(v).~\blacksquare$

\parindent=8mm From the above theorem we have the following corollaries.\\

\noindent
{\bf{Corollary 3.2}} {\it{A semi-slant warped product submanifold $M=N_T\times{_{f}N_\theta}$ of a nearly cosymplectic manifold is simply a Riemannian product of $N_T$ and $N_\theta$ if and only if $h(X, Z)\in\Gamma(\mu)$, for all  X$\in\Gamma(TN_T)$ and Z$\in\Gamma(TN_\theta)$, where $\mu$ is the invariant normal subbundle of $T^\perp M$.}}\\

\noindent
{\bf{Corollary 3.3}} {\it{There does not exist a mixed geodesic warped product semi-slant submanifold of a nearly cosymplectic manifold.}}

\parindent=8mm From Lemma 3.2 $(i),~(iv)$ and Theorem 3.4, we obtain
$$g(h(\xi, Z), FPZ)=g(h(\xi, PZ), FZ)=0,\eqno(3.9)$$
for any $Z\in\Gamma(TN_\theta)$.

\parindent=8mm Using the previous results, we derive the following geometric inequality for the length of second fundamental form.\\

{\bf{Theorem 3.5}} {\it Let M=$N_T\times {_{f} N_\theta}$ be a warped product semi-slant submanifold of a nearly cosymplectic manifold $\bar M$, where $N_T$ and $N_\theta$ are invariant and proper slant submanifolds of $\bar M$, respectively. Then}
\begin{enumerate}
\item[{(i)}] {\it{The squared norm of the second fundamental form $h$ satisfies}}
$$\|h\|^2 \geq 4q\csc^2\theta\{1+\frac{1}{9}\cos^2\theta\}\|\nabla \ln f\|^2\eqno(3.10)$$
{\it{where $\nabla \ln f$ is the gradient of $\ln f$ and $2q$ is the dimension of $N_\theta.$}}
\item[{(ii)}] {\it{If the equality sign of (3.10) holds, then $N_T$ is totally geodesic in $\bar M$ and $N_\theta$ is a totally umbilical submanifold of $\bar M$. Moreover, $M$ is a minimal submanifold of $\bar M$.}}
\end{enumerate}

\noindent
{\it {Proof.}} Let $\bar M$ be a $(2n+1)$-dimensional nearly cosymplectic manifold and let $M=N_T\times {_{f} N_\theta}$ be a warped product semi-slant submanifold of $\bar M$. Let us denote by ${\cal D}$ and ${\cal D}_\theta$ the tangent bundles on $N_T$ and $N_\theta$, respectively and let $dim N_T=2p+1$ and $dim N_\theta=2q$, then $m=2p+2q+1$. Let $\{X_1,\cdots, X_p, X_{p+1}=\phi X_1,\cdots, X_{2p}=\phi X_p, X_{2p+1}=\xi\}$ and $\{ Z_1,\cdots, Z_{q} , Z_{q+1}=\sec\theta PZ_1,\cdots, Z_{2q}=\sec\theta PZ_{q}\}$ be the local orthonormal frames of ${\cal D}$ and ${\cal D}_\theta$, respectively. Then, the orthonormal frames of $F{\cal D}_\theta$ and $\mu$ are $\{Z_1^\star=\csc\theta FZ_1,\cdots, Z_q^\star=\csc\theta FZ_{q}, Z_{q+1}^\star=\csc\theta\sec\theta FPZ_1,\cdots, Z_{2q}^\star=\csc\theta\sec\theta FPZ_{q}\}$ and $\{N_{m+2q+1}^\star\cdots, N_{2n+1}^\star\}$, respectively, where $N$ is the normal vector in the invariant normal subbundle of $T^\perp M$. The dimensions of $F{\cal D}_\theta$ and $\mu$ will be $2q$ and $2n+1-m-2q$, respectively. The length of second fundamental form $h$ is defined as
$$\|h\|^2=\sum_{r=m+1}^{2n+1}\sum_{i, j=1}^{m} g(h(U_i, V_j), N_r)^2\eqno(3.11)$$
for any vector fields $U_i, V_j$ tangent to $M$ and $N_r$ normal to $M$. Now, for the assumed frames of $F{\cal D}_\theta$ and $\mu$, the above equation can be written as
$$\|h\|^2=\sum_{r=m+1}^{m+2q}\sum_{i, j=1}^{m} g(h(U_i, V_j), Z_r^\star)^2+\sum_{r=m+2q+1}^{2n+1}\sum_{i, j=1}^{m}g(h(U_i, V_j), N_r^\star)^2.\eqno(3.12)$$
The first term in the right hand side of the above equation is the $F{\cal D}_\theta$-component and the second term is $\mu$-component. Here, we equate only $F{\cal D}_\theta$-component term, thus we have
$$\|h\|^2\geq\sum_{r=m+1}^{m+2q}\sum_{i, j=1}^{m} g(h(U_i, V_j), Z_r^\star)^2.\eqno(3.13)$$
Thus, for the given frames of ${\cal D}$ and ${\cal D}_\theta$, the above equation will be
$$\|h\|^2\geq\sum_{r=1}^{2q} \sum_{i, j=1}^{2p+1} g(h(X_i, X_j), Z_r^\star)^2+2\sum_{r=1}^{2q}\sum_{i=1}^{2p+1}\sum_{ j=1}^{2q} g(h(X_i, Z_j), Z_r^\star)^2$$
$$~~~~+\sum_{r=1}^{2q} \sum_{i, j=1}^{2q} g(h(Z_i, Z_j), Z_r^\star)^2.\eqno(3.14)$$
By Lemma 3.2 (ii), the first term of the right hand side is identically zero and we shall compute the next term and leave the third term
$$\|h\|^2\geq2\sum_{r=1}^{2q}\sum_{i=1}^{2p+1}\sum_{j=1}^{2q}g(h(X_i, Z_j), Z_r^\star)^2.$$
As $j, r=1,\cdots, 2q$, then the above equation can be written for one summation as
$$\|h\|^2\geq2\sum_{i=1}^{2p+1}\sum_{j=1}^{2q}g(h(X_i, Z_j), Z_j^\star)^2.$$
Separating the $h(\xi, Z)$-components, the above inequality will be
$$\|h\|^2\geq2\sum_{i=1}^{2p}\sum_{j=1}^{2q}g(h(X_i, Z_j), Z_j^\star)^2+2\sum_{j=1}^{2q}g(h(\xi, Z_j), Z_j^\star)^2.\eqno(3.15)$$
Now, we solve the second term of right hand side of (3.15) as follows
$$\sum_{j=1}^{2q}g(h(\xi, Z_j), Z_j^\star)^2=\csc^2\theta\sum_{j=1}^{q}g(h(\xi, Z_j), FZ_j)^2$$
$$~~~~~~~~~~~~+\csc^2\theta\sec^4\theta\sum_{j=1}^{q}g(h(\xi, PZ_j), FPZ_j)^2$$
$$~~~~~~~~~~+\csc^2\theta\sec^2\theta\sum_{j=1}^{q}g(h(\xi, Z_j), FPZ_j)^2$$
$$~~~~~~~~~~~+\csc^2\theta\sec^2\theta\sum_{j=1}^{q}g(h(\xi, PZ_j), FZ_j)^2.\eqno(3.16)$$
From (3.9), the last two terms of right hand side of (3.16) are identically zero and we will compute the first two terms as follows. We know that
$$g(h(\xi, Z), FZ)=g(\bar\nabla_Z\xi, FZ)=-g(\xi, \bar\nabla_ZFZ).$$
Using (2.7) and then (2.9), we obtain
$$g(h(\xi, Z), FZ)=-g(\xi, (\bar\nabla_Z\phi)Z)-g(\xi, \phi\bar\nabla_ZZ)+g(\xi, \bar\nabla_ZPZ).$$
Using the nearly cosymplectic character the first term of right hand side is zero, second and last terms are also zero by using (2.1), property of Riemannian connection and either orthogonality of vectors $Z$ and $PZ$ or the fact that $\xi\ln f=0$ or both, hence
$$g(h(\xi, Z), FZ)=0.\eqno(3.17)$$
If we interchange $Z$ by $PZ$, then $g(h(\xi, PZ), FPZ)=0.$ Put all these values in (3.16), we obtain
$$\sum_{j=1}^{2q}g(h(\xi, Z_j), FZ_j)^2=0.\eqno(3.18)$$
Thus, from (3.15) and (3.18), we derive
$$\|h\|^2\geq2\sum_{i=1}^{2p}\sum_{j=1}^{2q}g(h(X_i, Z_j), Z_j^\star)^2.$$
Using the frame of $F{\cal D}_\theta$, the above inequality can be written as
$$\|h\|^2\geq2\csc^2\theta\sum_{i=1}^{p}\sum_{j=1}^{q}g(h(X_i, Z_j), FZ_j)^2~~~~~~~~~~$$
$$~~~~~~~~~~~+2\csc^2\theta\sec^4\theta\sum_{i=1}^{p}\sum_{j=1}^{q}g(h(X_i, PZ_j), FPZ_j)^2$$
$$+2\csc^2\theta\sum_{i=1}^{p}\sum_{j=1}^{q}g(h(\phi X_i, Z_j), FZ_j)^2~$$
$$~~~~~~~~~~~~~~+2\csc^2\theta\sec^4\theta\sum_{i=1}^{p}\sum_{j=1}^{q}g(h(\phi X_i, PZ_j), FPZ_j)^2~~$$
$$~~~~~~~~~~~+2\csc^2\theta\sec^2\theta\sum_{i=1}^{p}\sum_{j=1}^{q}g(h(X_i, PZ_j), FZ_j)^2$$
$$~~~~~~~~~~~+2\csc^2\theta\sec^2\theta\sum_{i=1}^{p}\sum_{j=1}^{q}g(h(X_i, Z_j), FPZ_j)^2~~~~~~~$$
$$~~~~~~~~~~~+2\csc^2\theta\sec^2\theta\sum_{i=1}^{p}\sum_{j=1}^{q}g(h(\phi X_i, PZ_j), FZ_j)^2$$
$$~~~~~~~~~~~+2\csc^2\theta\sec^2\theta\sum_{i=1}^{p}\sum_{j=1}^{q}g(h(\phi X_i, Z_j), FPZ_j)^2.\eqno(3.19)$$
The first four terms of above inequality will be solved as follows. From Lemma 3.2 (iii), we have
$$g(h(\phi X, Z), FZ)=(X\ln f)g(Z, Z).\eqno(3.20)$$
Interchanging $X$ by $\phi X$ and then using (2.1), we derive
$$-g(h(X, Z), FZ)+\eta(X)g(h(\xi, Z), FZ)=(\phi X\ln f)g(Z, Z).$$
But from (3.17), the second term of left hand side of above equation vanishes identically, thus we obtain
$$g(h(X, Z), FZ)=-(\phi X\ln f)g(Z, Z).\eqno(3.21)$$
Interchanging $Z$ by $PZ$ in (3.20) and (3.21) and using (2.21) and the fact that $\xi$ is tangent to $N_T$, thus we obtain the following equations, respectively
$$g(h(\phi X, PZ), FPZ)=(X\ln f)\cos^2\theta g(Z, Z)\eqno(3.22)$$
and
$$g(h(X, PZ), FPZ)=-(\phi X\ln f)\cos^2\theta g(Z, Z).\eqno(3.23)$$
The last four terms of (3.19) will be solved by Theorem 3.4 as follows
$$g(h(X, Z), FPZ)=-g(h(X, PZ), FZ)=\frac{1}{3}(X\ln f)\cos^2\theta g(Z, Z).\eqno(3.24)$$
Interchanging $X$ by $\phi X$ in (3.24), we obtain
$$g(h(\phi X, Z), FPZ)=-g(h(\phi X, PZ), FZ)=\frac{1}{3}(\phi X\ln f)\cos^2\theta g(Z, Z).\eqno(3.25)$$
Put all these values in (3.19), we derive
$$\|h\|^2\geq4\csc^2\theta\sum_{i=1}^{p}\sum_{j=1}^{q}\{1+\frac{1}{9}\cos^2\theta\}(X_i\ln f)^2g(Z_j, Z_j)^2$$
$$~~~~~~~~+4\csc^2\theta\sum_{i=1}^{p}\sum_{j=1}^{q}\{1+\frac{1}{9}\cos^2\theta\}(\phi X_i\ln f)^2g(Z_j, Z_j)^2.$$
Hence, from (3.3) the above expression will be
$$\|h\|^2\geq4\csc^2\theta\sum_{j=1}^{q}\{1+\frac{1}{9}\cos^2\theta\}\|\nabla\ln f\|^2g(Z_j, Z_j)^2$$
$$~~~~~~~~~~~=4q\csc^2\theta\{1+\frac{1}{9}\cos^2\theta\}\|\nabla\ln f\|^2,$$
which is the inequality (3.10). If the equality holds in (3.10), then by (3.12) and (3.14), we obtain
$$h({\cal D},{\cal D})=0,~~~h({\cal D}_\theta, {\cal D}_\theta)=0,~~h({\cal D}, {\cal D}_\theta)\subset F{\cal D}_\theta.\eqno(3.26)$$
Now, for any $Z, W\in\Gamma({\cal D}_\theta)$ and $X\in\Gamma({\cal D})$, we have
$$g(h^\star(Z, W), X)=g(\nabla_ZW, X)=-g(W, \nabla_ZX)=-(X\ln f)g(Z, W)$$
where $h^\star$ is the second fundamental form of $N_\theta$ in $M$. Using (3.2), we obtain
$$g(h^\star(Z, W), X)=-g(Z, W)g(\nabla\ln f, X),$$
where $\nabla\ln f$ is the gradient of $\ln f$. Thus from the last relation, we get
$$h^\star(Z, W)=-\nabla\ln fg(Z, W).\eqno(3.27)$$
Since $N_T$ is totally geodesic submanifold in $M$ (by Lemma 3.1 (i)), using this fact with the first condition of (3.26), we get $N_T$ is totally geodesic in $\bar M$. Also, the second condition of (3.26) with (3.27) implies that $N_\theta$ is totally umbilical in $\bar M$. Moreover all conditions of (3.26) imply that $M$ is a minimal submanifold of $\bar M$. This completes the proof of the theorem.~$\blacksquare$\\

\parindent=8mm If we consider $\theta=\frac{\pi}{2}$, then the inequality (3.10) generalizes the inequality which we have obtained for contact CR-warped products in [16].

\bigskip

Author' address:\\

Siraj Uddin

\noindent Institute of Mathematical Sciences, Faculty of Science,University of Malaya, 50603 Kuala Lumpur, Malaysia

\noindent {\it E-mail}: {\tt siraj.ch@gmail.com}\\

\medskip

Abdulqader Mustafa

\noindent Institute of Mathematical Sciences, Faculty of Science,University of Malaya, 50603 Kuala Lumpur, Malaysia

\noindent {\it E-mail}: {\tt abdulqader.mustafa@yahoo.com}\\

\medskip

Berardine Renaldo Wong

\noindent Institute of Mathematical Sciences, Faculty of Science,University of Malaya, 50603 Kuala Lumpur, Malaysia

\noindent {\it E-mail}: {\tt bernardr@um.edu.my}\\

\medskip

Cenap Ozel

\noindent Department of Mathematics, Abant Izzet Baysal University, 14280 Bolu, Turkey

{\it E-mail}: {\tt cenap.ozel@gmail.com}


\begin{thebibliography}{99}
\bibitem{A} M. Atceken, {\it{Contact CR-warped product submanifolds in cosymplectic space forms}}, Collect. Math. 62 (2011), 17-26.
\bibitem{B} R.L. Bishop and B. O'Neill, {\it Manifolds of negative curvature}, Trans. Amer. Math. Soc. 145 (1969), 1-49.
\bibitem{B} D.E. Blair, {\it{Almost contact manifolds with killing structure tensors I}}, Pac. J. Math. 39 (1971) 285-292.
\bibitem{B} D.E. Blair and K. Yano, {\it{Affine almost contact manifolds and f-manifolds with affine killing structure tensors}}, Kodai Math. Sem. Rep. 23 (1971) 473-479.
\bibitem{B} D.E. Blair and D.K. Showers, {\it{Almost contact manifolds with killing structures tensors II}}, J. Differ. Geom. 9 (1974), 577-582.
\bibitem{C} J.L. Cabrerizo, A. Carriazo, L.M. Fernandez and M. Fernandez, {\it Semi-slant submanifolds of a Sasakian manifold}, Geom. Dedicata 78 (1999), 183-199.
\bibitem{C} J.L. Cabrerizo, A. Carriazo, L.M. Fernandez and M. Fernandez, {\it{Slant submanifolds in Sasakian manifolds}}, Glasgow Math. J. 42 (2000), 125-138.
\bibitem{C} B.Y. Chen, {\it Geometry of warped product CR-submanifolds in Kaehler manifolds I}, Monatsh. Math. 133 (2001), 177-195.
\bibitem{H} I. Hasegawa and I. Mihai, {\it Contact CR-warped product submanifolds in Sasakian manifolds}, Geom. Dedicata 102 (2003), 143-150.
\bibitem{K} K.A. Khan, V.A. Khan and Siraj Uddin, {\it{Warped product submanifolds of cosymplectic manifolds}}, Balkan J. Geom. Appl. 13 (2008), 55-65.
\bibitem{B} B. O'Neill, {\it Semi-Riemannian geometry with application to Relativity}, Academic Press, New York, 1983.
\bibitem{P} N. Papaghiuc, {\it{Semi-slant submanifolds of Khaelerian manifolds}}, Ann. St. Univ. Iasi, 9 (1994), 55-61.
\bibitem{S} Siraj Uddin, V.A. Khan and K.A. Khan, {\it{A note on warped product submanifolds of cosymplectic manifolds}}, Filomat 24 (2011), 95-102.
\bibitem{S} Siraj Uddin and K.A. Khan, {\it{Warped product CR-submanifolds of cosymplectic manifolds}}, Ricerche di matematica 60 (2011), 143-149
\bibitem{S} Siraj Uddin, S. H. Kon, M. A. Khan and Khushwant singh, {\it{Warped product semi-invariant submanifolds of nearly cosymplectic manifolds}}, Math. Probl. Eng. 2011 (2011).
\bibitem{S} Siraj Uddin and Khalid Ali Khan, {\it{An inequality for contact CR-warped product submanifolds of nearly cosymplectic manifolds}}, J. Inequal. Appl. 304 (2012).


\end{thebibliography}
\end{document}